\documentclass[a4paper,11pt]{article}
\usepackage{geometry,amssymb,amsmath,enumerate,latexsym,array}
\usepackage[dvips]{graphicx}
\usepackage[all]{xy}
\newcommand{\id}{i\!d}
\newcommand{\ep}{\varepsilon}
\newcommand{\bA}{\mathsf{A}}
\newcommand{\bB}{\mathsf{B}}
\newcommand{\bC}{\mathsf{C}}

\newcommand{\bP}{\mathsf{P}}
\newcommand{\sets}{\mathsf{Sets}}

\newcommand{\lra}{\leftrightarrow}
\newcommand{\lan}{\langle}
\newcommand{\ran}{\rangle}

\def\ob{{\mathrm{ Ob}}}

\def\bar{\overline}

\def\leq{\leqslant}

\newcommand{\arr}{\mathrm{Arr}}
\newtheorem{example}{Example}[section]

\newtheorem{thm}[example]{Theorem}

\begin{document}

\begin{center}
\textbf{\LARGE Rewriting procedures generalise to\\
\vspace{2mm}
Kan extensions of actions of categories}
\end{center}

\section{Introduction}

This is a brief account of work of Brown and Heyworth
\cite{paper2} on extensions of rewriting methods.

The standard expression of such methods is in terms of words $w$ in a
free monoid $\Delta^*$ on a set $\Delta$. This may be extended to terms
$x|w$ where $x$ belongs to a set $X$ and the link between $x$ and $w$ is in
terms of an action. More precisely, we suppose a monoid $A$ acts on the
set $X$ on the right, and there is given a morphism of monoids $F$:
$A\to B$ where $B$ is given by a presentation with generating set
$\Delta$. The result of the rewriting will then be normal forms for
the {\it induced action} of $B$ on $F_*(X)$. This gives an important 
extension of rewrite methods. 

In fact monoids may be replaced by categories, and sets by directed
graphs. This gives a formulation in terms of Kan extensions, or
induced actions of categories, which  we now explain.

\section{Presentations of Kan Extensions}

Let $\bA$ be a category.
  A \textbf{category action} $X$ of $\bA$ is a 
functor $X:\bA \to\sets$.
Let $\bB$be a second category and let $F:\bA \to \bB$ be a  functor. 
Then an\textbf{extension of the action $X$ along $F$} is a pair $(K,\ep)$
where $K:\bB \to \sets$ is a  functor and $\ep:X \to K \circ F$ is
a natural transformation.
The \textbf{Kan extension of the action $X$ along $F$} is an
extension of the action $(K,\ep)$ with the universal property that
for any other extension of the action $(K',\ep')$ there exists a
unique natural transformation $\alpha:K \to K'$ such that
$\ep'=\alpha \circ \ep$.

The problem that has been introduced is that of ``computing a Kan
extension''. In order to keep the analogy with computation and
rewriting for presentations of monoids we propose a definition of
a presentation of a Kan extension. The papers 
\cite{BLW, CaWa1, CaWa2, Rosebrugh} were very influential on the 
current work.

A \textbf{Kan extension data} $(X',F')$ consists of small categories
$\bA$, $\bB$ and functors $X':\bA \to \sets$ and $F':\bA \to \bB$.
A \textbf{Kan extension presentation} is a quintuple
$\mathcal{P}:=kan \lan \Gamma | \Delta | RelB |X | F \ran$ where
\begin{enumerate}
  \item $\Gamma$ and $\Delta$ are (directed) graphs;
  \item $X: \Gamma \to \sets$ and $F: \Gamma \to P \Delta$ are graph
morphisms to the category of sets and the free category on
$\Delta$ respectively;
  \item  and $RelB$ is a set of relations on the free category $P\Delta$.
\end{enumerate}

We say $\mathcal{P}$ \textbf{presents} the Kan extension $(K,\ep)$ 
of the Kan extension
data $(X',F')$ where $X':\bA \to \sets$ and $F':\bA \to \bB$ if
\begin{enumerate}
\item
 $\Gamma$ is a generating 
graph for $\bA$  and $X: \Gamma \to \sets$ is the restriction of
 $X':\bA \to \sets$ 
\item
$cat \lan \Delta | RelB \ran$ is a category presentation for $\bB$.
\item
$F:\Gamma \to P\Delta$ induces $F':\bA \to \bB$.
\end{enumerate}

We expect that a Kan extension $(K, \ep)$ is given by a set $KB$
for each $B \in \ob \Delta$ and a function $Kb: KB_1 \to KB_2$ for
each $b:B_1 \to B_2 \in \bB$  (defining the functor $K$) together
with a function $\ep_A: XA \to KFA$ for each $A \in \ob \bA$ (the
natural transformation). This information can be given in four
parts:
\begin{itemize}
  \item  the set $\bigsqcup_B  KB$;
  \item  a function $\bar{\tau}:\bigsqcup_B KB \to \ob \bB$;
  \item  a partial function (action) $\bigsqcup_B KB \times \arr\bP \to \bigsqcup_B
  KB$;
  \item  and a function $\ep: \bigsqcup_A XA \to \bigsqcup_B  KB$.
\end{itemize}
Here $\bigsqcup_B  KB$ and $\bigsqcup_A XA$ are the disjoint
unions of the sets $KB$, $XA$ over $\ob \bB$, $\ob \bA$
respectively; if $z \in KB$ then $\bar{\tau}(z)=B$
 and if further $src(p)=B$ for $p \in \arr\bP$ then $z \cdot p$
is defined.

\section{Rewriting for Kan Extensions}

 The main result of the paper defines rewriting procedures on the 
$P \Delta$-set 
$$
T:= \bigsqcup_{B \in \ob\Delta} \bigsqcup_{A \in
\ob \Gamma} XA \times P\Delta(FA,B). 
$$ 
Two kinds of rewriting are involved here. The
first is the familiar $x|ulv \to x|urv$ given by a relation
$(l,r)$. The second derives from a given action of certain words
on elements, so allowing rewriting $x|F(a)v \to x \cdot a|v$. 
Further,the elements $x$ and $x \cdot a$ may belong to different 
sets. When such rewritingprocedures complete, the associated normal form 
gives in effect acomputation of what we call the 
{\em Kan extension defined by the presentation}.

\begin{thm}
Let $\mathcal{P}=kan \lan \Gamma | \Delta | RelB | X F \ran$ be a
Kan extension presentation, and let $\bP$, $T$, $R=(R_\ep, R_K)$
be defined as above. Then the Kan extension $(K,\ep)$ presented by
$\mathcal{P}$ may be  given by the following data:
\begin{enumerate}[i)]
\item
the set $\bigsqcup_B KB = T/ \stackrel{*}{\lra}_R$,
\item
the function $\bar{\tau}: \bigsqcup_B KB \to \ob \bB$ induced by
$\tau:T \to \ob \bP$,
\item
the action of $\bB$ on $ \bigsqcup_B KB$ induced by the action of
$\bP$ on $T$,
\item
the natural transformation $\ep$ determined by $x \mapsto [x|
\id_{FA}]$ for $x \in XA, \;A \in \ob \bA$.
\end{enumerate}
\end{thm}

\subsection{Reduction and critical pairs}
To work with a rewrite system $R$ on $T$ we will require certain
concepts of order on $T$. We  give properties of  orderings $>_X$
on $\bigsqcup_A XA$ and  $>_P$ on $\arr\bP$ to enable us to
construct an ordering $>_T$ on $T$ with the properties needed for
the rewriting procedures.

Given an admissible well-ordering $>_T$ on $T$ it is
possible to discuss when a reduction relation generated by a
rewrite system is compatible with this ordering.
It is a standard result that if a reduction relation is compatible
with an admissible well-ordering, then it is Noetherian. 
A Noetherian reduction relation on a set is confluent if it is
locally confluent(Newman's Lemma\cite{TAT}).
By standard abuse of
notation the rewrite system $R$ will be called complete when $\to_R$ is
complete. 
Hence, if $R$ is compatible with an admissible well-ordering on
$T$ and  $\to_R$ is locally confluent then $\to_R$ is complete. By
orienting the pairs of $R$ with respect to the chosen ordering
$>_T$ on $T$, $R$ is made to be Noetherian. The problem remaining
is testing for local confluence of $\to_R$ and changing $R$ in
order to obtain an equivalent confluent reduction relation.

We explain the notion of critical pair for a rewrite
system for $T$, extending the traditional notion to our situation.
In particular the overlaps involve either just $R_T$, or just
$R_P$ or an interaction between $R_T$ and $R_P$.

A term $crit\in T$ is called \textbf{critical} if it may be reduced by
two or more different rules.
A pair $(crit1,crit2)$ of distinct terms resulting from two single-step
reductions of the same term is called a \textbf{critical pair}.
A critical pair for a reduction relation $\to_R$ is said to
\textbf{resolve} if
there exists a term $res$ such that both $crit1$ and $crit2$ reduce to a
$res$, i.e. $crit1 \stackrel{*}{\to}_R res$,
$crit2 \stackrel{*}{\to}_R res$.
If $t=x|b_1 \cdots b_n$, then a \textbf{part} of $t$ is either a term 
$x|b_1 \cdots b_i$ for some $1 \leq i \leq n$ or a word $b_i b_{i+1} \cdots b_j$
for some $1 \leq i \leq j \leq n$.
Let $(rule1,rule2)$ be a pair of rules of the rewrite
system $R=(R_T,R_P)$ where $R_T \subseteq T \times T$ and $R_P
\subseteq \arr\bP \times \arr \bP$. Then the rules are said to \textbf{overlap}
when $rule1$ and $rule2$ may both be applied to the same term $t$ in such a way
that there is a part $c$ of the term that is affected by both the rules.

There are five types of overlap for this kind of rewrite system,
as shown in the following table:
\vspace{2ex}
\begin{center}
\setlength{\extrarowheight}{6pt}
\begin{tabular}{||c|c|c|c|c|l|l||}
\hline
 \# &rule1 & in & rule2& in  & \hspace{1cm} overlap & critical pair \\
\hline
(i) &$(s_1,u_1)$ & $ R_T$  & $(s_2,u_2)$& $ R_T$  &  $s_2=s_1\cdot q$ for some $q \in \arr \bP$ & $(u_1\cdot q,u_2 )$ \\
\hline
(ii)&$(l_1, r_1)$ & $R_P$ & $(l_2,r_2)$ & $ R_P$   &  $l_1=pl_2q$ for some $p,q \in \arr \bP$    &  $(r_1,pr_2q)$
\\ \cline{1-1} \cline{6-7}
(iii) & & &  &   &  $l_1q=pl_2$  for some $p,q \in \arr \bP$   &   $(r_1q,pr_2)$
\\ \hline
(iv)& $(s_1,u_1)$ & $ R_T$ & $(l_1, r_1)$ & $R_P$  &  $s_1\cdot q= s\cdot l_1$ for some $s \in T, q \in \arr \bP$  & $(u_1\cdot q, s \cdot r_1)$
 \\\cline{1-1}  \cline{6-7}
 (v)&  & & &  &   $s_1 =s \cdot (l_1q)$ for some $s \in T, q \in \arr \bP$      & $(u_1, s \cdot r_1q)$  \\
  \hline
\end{tabular}

Overlap table

\end{center}

\subsection{Completion procedure}
We show:
(i) how to find overlaps between rules of $R$;
(ii) how to test whetherthe resulting critical pairs resolve;
(iii) that if all the critical pairs resolve
      then this implies $\to_R$ is confluent; and 
(iv) that critical pairswhich do not resolve may be added to $R$ 
without affecting theequivalence relation $R$ defines on $T$.
We have now set up and proved everything necessary for a variant of the
Knuth-Bendix procedure, which will add rules to a rewrite system $R$
resulting from a presentation of a Kan extension, to attempt to find an
equivalent complete rewrite system $R^C$. The benefit of such a system is that
$\to_{R^C}$ then acts as a normal form function for $\stackrel{*}{\lra}_{R^C}$
on $T$.

\begin{thm}
\label{proc} Let $\mathcal{P}= \lan \Gamma | \Delta | RelB |X |F
\ran$ be a finite presentation of a Kan extension $(K, \ep)$. Let
$P:=P \Delta$, $T:= \bigsqcup_{\ob\Delta} \bigsqcup_{\ob\Gamma} XA \times
\bP(FA, B)$, and let $R$ be the initial rewrite system for
$\mathcal{P}$ on $T$. Let $>_T$ be an admissible well-ordering on
$T$. Then there exists a procedure which, if it terminates,  will
return a rewrite system $R^C$ which is complete with respect to
the ordering $>_T$ and  such that the equivalence relations
$\stackrel{*}{\lra}_R$, $\stackrel{*}{\lra}_{R^C}$ coincide.
\end{thm}

The above  procedure which attempts completion of a presentation of a
Kan extension has been implemented in {\sf GAP3}.

\section{Example of a {\sf GAP} Session on the Rewriting Procedure}
\label{kaneg}

Here we give an example to show the use of the implementation.
Let $\bA$ and $\bB$ be the categories generated by the graphs
below, where $\bB$ has the relation $b_1b_2b_3=b_4$. {\large
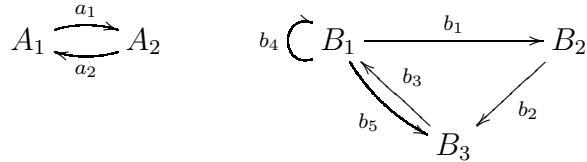
$$
\xymatrix{A_1 \ar@/^/[r]^{a_1}
        & A_2 \ar@/^/[l]^{a_2}
       && B_1  \ar@(dl,ul)^{b_4}
             \ar[rr]^{b_1}
             \ar@/_/[dr]_{b_5}
       && B_2 \ar[dl]^{b_2} \\
     &&&& B_3 \ar[ul]_{b_3} \\}
$$}
Let $X:\bA \to \sets$ be defined by $ XA_1 = \{ x_1, x_2, x_3 \},
                                   \ XA_2 = \{ y_1, y_2 \}$   with\\
$Xa_1:XA_1 \to XA_2: x_1 \mapsto y_1, x_2 \mapsto y_2, x_3 \mapsto
y_1$,\\ $Xa_2:XA_1 \to XA_2: y_1 \mapsto x_1, y_2 \mapsto x_2,$\\
and let $F:\bA \to \bB$ be defined by $FA_1=B_1, \ FA_2=B_2, \
Fa_1=b_1$ and $Fa_2 = b_3 b_2$. The input to the computer program
takes the following form.
First read in the program and set up the variables:
\begin{verbatim}
gap> RequirePackage("kan");
gap> F:=FreeGroup("b1","b2","b3","b4","b5","x1","x2","x3","y1","y2");;
gap> b1:=F.1;;b2:=F.2;;b3:=F.3;;b4:=F.4;;b5:=F.5;;
gap> x1:=F.6;;x2:=F.7;;x3:=F.8;;y1:=F.9;;y2:=F.10;;
\end{verbatim}
\newpage
Then we input the data (choice of names is unimportant):
\begin{verbatim}
gap> OBJa:=[1,2];;
gap> ARRa:=[[1,2],[2,1]];;
gap> OBJb:=[1,2,3];;
gap> ARRb:=[[b1,1,2],[b2,2,3],[b3,3,1],[b4,1,1],[b5,1,3]];;
gap> RELb:=[[b1*b2*b3,b4]];;
gap> fOBa:=[1,2];;
gap> fARRa:=[b1,b2*b3];;
gap> xOBa:=[[x1,x2,x3],[y1,y2]];;
gap> xARRa:=[[y1,y2,y1],[x1,x2]];;
\end{verbatim}
To combine all this data in one record (the field names are important):
\begin{verbatim}
gap> KAN:=rec( ObA:=OBJa, ArrA:=ARRa,  ObB:=OBJb, ArrB:=ARRb, RelB:=RELb,
              FObA:=fOBa, FArrA:=fARRa, XObA:=xOBa, XArrA:=xARRa );;
\end{verbatim}
To calculate the initial rules do:
\begin{verbatim}
gap> InitialRules( KAN );
\end{verbatim}
The output will be:
\begin{verbatim}
i= 1, XA= [ x1, x2, x3 ], Ax= x1, rule= [ x1*b1, y1 ]
i= 1, XA= [ x1, x2, x3 ], Ax= x2, rule= [ x2*b1, y2 ]
i= 1, XA= [ x1, x2, x3 ], Ax= x3, rule= [ x3*b1, y1 ]
i= 2, XA= [ y1, y2 ], Ax= y1, rule= [ y1*b2*b3, x1 ]
i= 2, XA= [ y1, y2 ], Ax= y2, rule= [ y2*b2*b3, x2 ]
[ [ b1*b2*b3, b4 ], [ x1*b1, y1 ], [ x2*b1, y2 ], [ x3*b1, y1 ],
  [ y1*b2*b3, x1 ], [ y2*b2*b3, x2 ] ]
\end{verbatim}
This means that there are five initial $\ep$-rules:\\
$( \ x_1|Fa_1,  x_1.a_1|\id_{FA_2} \ ), \
 ( \ x_2|Fa_1,  x_2.a_1|\id_{FA_2} \ ),
( \ x_3|Fa_1,  x_3.a_1|\id_{FA_2} \ ),$\\ $
 ( \ y_1|Fa_2,  y_1.a_1|\id_{FA_1} \ ), \
 ( \ y_2|Fa_2,  y_2.|a_1 \id_{FA_1} \ ),$ \\
i.e. $ \ x_1|b_1 \to y_1|\id_{B_2}, \ x_2|b_1 \to y_2|\id_{B_2}, \
x_3|b_1 \to y_1|\id_{B_2}, \  y_1|b_2b_3  \to  x_1|\id_{B_1},
y_2|b_2b_3 \to x_2|\id_{B_1} \ $

and one initial $K$-rule: $b_1b_2b_3 \to  b_4$.

To attempt to complete the Kan extension presentation do:
\begin{verbatim}
gap> KB( InitialRules(KAN) );
\end{verbatim}
The output is:
\begin{verbatim}
[ [ x1*b1, y1 ], [ x1*b4, x1 ], [ x2*b1, y2 ], [ x2*b4, x2 ], [ x3*b1, y1 ],
  [ x3*b4, x1 ], [ b1*b2*b3, b4 ], [ y1*b2*b3, x1 ], [ y2*b2*b3, x2 ] ]
\end{verbatim}
In other words to complete the system we have to add the rules
$$x_1|b_4  \to  x_1, \quad x_2|b_4  \to  x_2, \text{ and } x_3|b_4  \to
x_1.$$
The result of attempting to compute the sets by doing:
\begin{verbatim}
gap> Kan(KAN);
\end{verbatim}
is a long list and then:
\begin{verbatim}
enumeration limit exceeded: complete rewrite system is:
[ [ x1*b1, y1 ], [ x1*b4, x1 ], [ x2*b1, y2 ], [ x2*b4, x2 ], [ x3*b1, y1 ],
  [ x3*b4, x1 ], [ b1*b2*b3, b4 ], [ y1*b2*b3, x1 ], [ y2*b2*b3, x2 ] ]
\end{verbatim}
This means that the sets $KB$ for $B$ in $\bB$ are too large. The
limit set in the program is  1000. (To change this the user should type
{\tt EnumerationLimit:= 5000} -- or whatever, after reading in the program.)
In fact the above example is
infinite.
The complete rewrite system is output instead of the
sets. We can in fact use this to obtain regular expressions for
the sets. In this case the regular expressions are:
\begin{center}
\begin{tabular}{lcl}
$KB_1$ & $:=$ & $(x_1+x_2+x_3)|(b_5(b_3{b_4}^*b_5)^*b_3{b_4}^*+\id_{B_1}).$\\
$KB_2$ & $:=$ & $(x_1+x_2+x_3)|b_5(b_3{b_4}^*b_5)^*b_3{b_4}^*(b_1) +
(y_1+y_2)|\id_{B_2}.$\\
$KB_3$ & $:=$ & $(x_1+x_2+x_3)|b_5(b_3{b_4}^*b_5)^*(b_3{b_4}^*b_1b_2
+\id_{B_3}) + (y_1 + y_2)|b_2.$\\
\end{tabular}
\end{center}
The actions of the arrows are defined by concatenation followed by
reduction.\\ For example $x_1|b_5b_3b_4b_4b_5$ is an element of
$KB_3$, so $b_3$ acts on it to give $x_1|b_5b_3b_4b_4b_5b_3$ which
is irreducible, and an element of $KB_1$.

The general method of  obtaining regular expressions for these
computations will be given in a separate paper (see Chapter 4 of
\cite{Anne}).

\section{Applications}
Mac\,Lane
wrote that ``the notion of Kan extensions subsumes all the other
fundamental concepts of category theory'' in section 10.7 of
\cite{Mac} (entitled ``All Concepts are Kan Extensions'').
So the power of rewriting theory may now be brought to bear on a
much wider range of combinatorial enumeration problems.
Traditionally rewriting is used for solving the word problem for
monoids. It has also been used for coset enumeration problems
\cite{Redfern,Hurt}.  It may now also be used in the specification of
\begin{enumerate}[i)]
\item equivalence classes and equivariant equivalence classes,
\item arrows of a category or groupoid,
\item right congruence classes given by a relation on a monoid,
\item orbits of an  action of a group or monoid.
\item conjugacy classes of a group,
\item coequalisers, pushouts and colimits of sets,
\item induced permutation representations of a group or monoid.
\end{enumerate}
and many others.
In this paper we are concerned with the description of the theory and
the implementation of the procedure with respect to one ordering.
It is hoped to consider implementation of efficiency strategies and other 
orderings onanother occasion.The advantages of our abstraction should 
then become even clearer, since one efficient implementation will be able to
apply to a variety of situations, including some not yet apparent.

\section{Further work, questions}

\subsection{Iteration}
One of the pleasant features of the procedure we describe is that the input and 
the output are of a similar form. The consequence of this is that if the action
$K$, given by $(X',F')$, has been defined on $\Delta$,
then given a second functor $G':\bB \to \bC$ and a presentation 
$cat\lan \Lambda| RelC\ran$ for $\bC$, it is straightforward to consider a 
presentation for the Kan extension data $(K,G')$. This new extension is in fact
the Kan extension with data $(X',G' \circ F')$.

\subsection{Kan Extensions and Noncommutative Gr\"obner Bases}
It is well-known that rewrite systems are a special case of noncommutative
Gr\"obner bases. It is possible to express a $K$-algebra presentation as an example
of a Kan extension over $K$-categories but it is not clear how to apply
Gr\"obner basis procedures to general Kan extensions of actions of 
$K$-categories.

\subsection{Orderings on $\bP$-sets}
In our paper we put stronger conditions on the ordering than may be necessary.
Weaker conditions may or may not have an advantage. The only ordering we have 
implemented is the standard length-lexicographical. The choice of orderings may
be wider than with ordinary rewriting, and this has not been investigated.

\subsection{Language Theory}
The actions in question are category actions on sets. Thus for each object of
$\bB$ we wish to specify a set, and for each arrow of $\bB$ we require a function
defined on the sets. In theory it is fine to specify sets by equivalence classes
of a larger set, with a normal form function. In practice we may wish to get 
hold of an expression for all the normal forms. When the sets are finite we can 
use a basic enumeration procedure, but when the sets are infinite, enumeration is 
not an answer. In this case an automaton can be constructed from the complete
rewrite system and language equations can be obtained and manipulated to obtain
a {\em regular expression} for the normal forms of the elements of each set.
It would be nice to program this!

\subsection{Automatic Kan Extensions}
Given the existing and current work on automatic groups, semigroups and coset 
systems it is natural to ask: what does the concept of automatic mean in terms
of a Kan extension?
An automatic coset system consists of ``a finite state automaton that provides
a name for each coset, and a set of finite state automata that allow these 
cosets to be multiplied by the group generators'' \cite{Redfern}.
We would expect therefore that an automatic Kan extension system would consist 
of a finite state automata for each set $KB$ that provides a name for each 
element of the set, and a finite state automaton for each arrow on $\Delta$ 
that allows the sets to be acted upon by the arrows of $\bB$.

\end{document}